\documentclass[graybox]{svmult}

\usepackage{mathptmx}       
\usepackage{helvet}         
\usepackage{courier}        
\usepackage{type1cm}        
\usepackage{makeidx}         
\usepackage{graphics}        
\usepackage{multicol}        
\usepackage[bottom]{footmisc}

\usepackage{etex}
\reserveinserts{30}
\usepackage{amsmath, amssymb}
\usepackage{tikz}
\usepackage{xcolor,calc}
\usetikzlibrary{matrix,arrows,decorations.pathmorphing}
\usepackage{amsmath,amscd}%
\usepackage{amsfonts}%
\usepackage{amssymb}%
\usepackage{txfonts}
\usepackage{color}
\usepackage{subfigure}
\usepackage{pgfplots}
\usepackage{cancel}
\usepackage{polynom}

\usepackage{cite}

\makeindex             

\newcommand{\kdifform}[2]{#1^{(#2)}} 
\newcommand{\incidencederivative}[2]{\mathsf{E}^{(#1,#2)}} 
\newcommand{\innerspace}[3]{\left(#1,#2\right)_{#3}} 
\newcommand{\matrixoperator}[1]{#1}
\newcommand{\inner}[2]{\left( #1, #2\right)} 

\newcommand{\manifold}[1]{\mathcal{#1}} 

\renewcommand{\div}[1]{\nabla \cdot #1} 

\newcommand{\figref}[1]{Figure~\ref{#1}} 

\newcommand{\realNumber}{\mathbb{R}}

\begin{document}

\title*{Mixed Mimetic Spectral Element method applied to Darcy's problem}
\author{Pedro Pinto Rebelo, Artur Palha and Marc Gerritsma}
\institute{ \email{\{P.J.PintoRebelo, A.Palha, M.I.Gerritsma\}@tudelft.nl} \at Delft University of Technology, Faculty of Aerospace Engineering, Aerodynamics Group, Delft,the Netherlands}
%
%
\maketitle

\abstract*{We present a discretization for Darcy's problem using the recently developed \textit{Mimetic Spectral Element Method} \cite{kreeft2011mimetic}. The gist lies in the exact discrete representation of integral relations. In this paper, an anisotropic flow through a porous medium is considered and a discretization of a full permeability tensor is presented. The performance of the method is evaluated on standard test problems,  converging at the same rate as the best possible approximation.}

\abstract{We present a discretization for Darcy's problem using the recently developed \textit{Mimetic Spectral Element Method} \cite{kreeft2011mimetic}. The gist lies in the exact discrete representation of integral relations. In this paper, an anisotropic flow through a porous medium is considered and a discretization of a full permeability tensor is presented. The performance of the method is evaluated on standard test problems,  converging at the same rate as the best possible approximation.}

\section{DARCY FLOW}
\label{Section::Introduction}

Anisotropic heterogeneous diffusion problems are ubiquitous across different scientific fields, such as, hydrogeology, oil reservoir simulation, plasma physics, biology, etc \cite{herbin2008benchmark}. Darcy's equation describes a steady pressure-driven flow through a porous medium where fluxes and pressure are linearly related,
\[
\text{div} \ \frac{\mathbb{K}}{\mu} \ \text{grad} \ p = \phi \overset{\mu = 1}{\longrightarrow}
\left \{ \begin{array}{llr}
\vec{u} - \text{grad} \ p = 0   & \text{in } \Omega \quad \quad \quad \quad \quad \quad & (1\mathrm{a})\\
\text{div} \ \vec{q} = \phi  & \text{in } \Omega & (1\mathrm{b}) \\
\vec{q} = \mathbb{K} \vec{u} & \text{in } \Omega & (1\mathrm{c})\\
\vec{q} = \vec{q_{0}}  & \text{in } \partial \Omega & (1\mathrm{d})\\
\end{array} \right .
\label{equation::darcyProblem}
\]
\setcounter{equation}{1}
where $\vec{u}$ is the fluid velocity, $p$ the pressure, $\vec{q}$ the mass flux and $\phi$ the prescribed source term. Without loss of generality let the viscosity, $\mu = 1$, and consider a permeability symmetric, positive definite tensor denoted by $\mathbb{K}$.

\begin{figure}[ht]
 \centering
\includegraphics[width=0.4\textwidth]{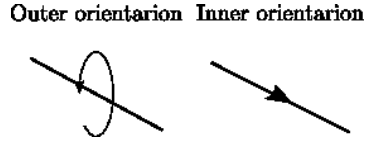}
\caption{Consider a line where we can have two types of orientation: \textit{Outer} - around the line; \textit{Inner} - along the line.}
\label{fig::lineOrientation}
\end{figure}

In a three-dimensional setting are: four types of submanifolds (\textit{points, lines, surfaces and volumes}); and two orientations (\textit{outer and inner}, as an example see \figref{fig::lineOrientation}). Tessellation divides the physical domain in a set of these geometric objects to which we associate discrete variables, i.e. \textit{integral quantities}. Thus, associated with every physical variable is a correspondent geometric object, this symbiotic relation between physics and geometry is the core of \textit{mimetic methods}. Many scholars are aware of this relationship \cite{bossavit2005discretization,burke1985applied,frankel,tonti1975formal}.

Starting from the mass balance equation, $(1\mathrm{b})$,
\begin{align}
 \int_{V} \text{div} \ \vec{q} \ \mathrm{d}V = \int_{\partial V} \vec{q} \cdot \vec{n} \ \mathrm{d}S = \int_{V} \phi \ \mathrm{d}V,
 \label{integral_divergence}
\end{align}
it is clear that the \textit{divergence in a volume} is equal to the sum of the \textit{surface integral quantities}, i.e. \textit{oriented fluxes}. Thus, we will associate mass fluxes, $\vec{q}$, with quantities that go \textit{through} surfaces. This equation therefore tells us that the right hand side term $\phi$ is associated to \textit{outer-oriented volumes}.

Similarly, using Newton-Leibniz relation for equation $(1\mathrm{a})$,
\begin{align}
 \int_{C} \text{grad} \ p \ \mathrm{d}C = \int_{\partial C} p = p\left(B \right) - p\left(A \right) = \int_{C} \vec{u} \ \mathrm{d}C,
\label{eq:u_grad_p}
\end{align}
the fluid velocity, $\vec{u}$, is represented {\em along} lines and $p$ is represented by the values in points. From (\ref{eq:u_grad_p}) we deduces that $u$ and $p$ are {\em inner-oriented variables}.

The \textit{constitutive/material relation} relation $(1\mathrm{c})$ is given by,
\begin{align}
\vec{q} = \mathbb{K} \vec{u},
\end{align}
which defines how quantities associated to inner-oriented lines relate to quantities associated to outer-oriented surfaces. Whereas equation
$(2)$ and  $(3)$ can be exactly satisfied on a finite grid the constitutive equation $(4)$ needs to be approximated.

\begin{figure}[ht]
 \centering
\includegraphics[width=0.65\textwidth]{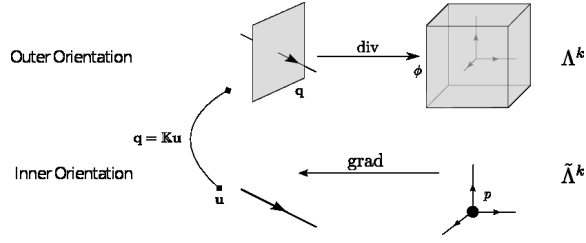}
\caption{Darcy's flow problem geometric characterization. Fluxes, $\vec{q}$, are associated with outer oriented surfaces; $\phi$, is associated with outer oriented volumes; velocity, $\vec{u}$, is associated with inner oriented lines; pressure, $p$, is associated with inner oriented points}
\label{fig::summaryDarcy}
\end{figure}

The importance of respecting the geometric nature in physics is discussed in \cite{gerritsmaICOSAHOM}. \figref{fig::summaryDarcy} summarizes the geometric character of the Darcy's problem.

We will denote the space of variables associated to outer-oriented $k$-dimensional objects by $\Lambda^k \left( \manifold{M} \right)$ and the space of variables associated to inner-oriented $k$-dimensional objects by $\tilde{\Lambda}^k \left( \manifold{M} \right)$ as indicated in Figure~\ref{fig::summaryDarcy}.

In this paper we will make use of the spectral element method described in \cite{gerritsmaICOSAHOM,kreeft2011mimetic}, application of these ideas to Stokes' flow see \cite{KreeftICOSAHOM,kreeft2012priori,kreeft::stokes}; Poisson equation for volume forms \cite{palhaLaplaceDualGrid}; advection equation \cite{palhaAdvection}; derivation of a momentum conservation scheme \cite{ToshniwalICOSAHOM}. Extension to compatible isogeometric methods see \cite{HiemstraICOSAHOM,hiemstra2012high}. For applications of these ideas in a finite difference setting see Brezzi et al. \cite{BrezziBuffaLipnikov2009}. In the context of finite element methods Arnold, Falk and Winther \cite{arnold2010finite} proposed a \textit{Finite Element Exterior Calculus}. In a more geometric spirit Desbrun et al. \cite{desbrun2005discrete} and Hirani \cite{Hirani_phd_2003} developed the \textit{discrete exterior calculus} (DEC). An application of the latter to Darcy flow can be found in \cite{hirani2008numerical}.

\section{DISCRETIZATION OF EQUATIONS}
\label{Section::Discretization}
In this section we will describe the discretization by defining the weak formulation. The approached followed here is similar to \cite{hyman1997numerical, jerome2012analysis}.

\begin{svgraybox}
For vectors associated with outer-oriented surfaces, $\Lambda^2\left(\manifold{M}\right )$, we define the weighted inner product
\begin{align}
 \innerspace{\mathbf{a}}{\mathbf{b}}{\manifold{M},\mathbb{K}} = \int_{\manifold{M}} \mathbf{a} \mathbb{K}^{-1}\mathbf{b} \ \mathrm{d}V.
\label{def:weighted_inner}
\end{align}
Furthermore, we define bilinear maps $((\cdot,\cdot))_{\manifold{M}}\,:\, \Lambda^1\left(\manifold{M}\right ) \times \tilde{\Lambda}^{n-1}\left(\manifold{M} \right ) \rightarrow \mathbb{R}$
\begin{align}
 ((\mathbf{u},{\mathbf{\tilde{v}}}))_{\manifold{M}} := \int_{\manifold{M}} \mathbf{u} \cdot \mathbf{\tilde{v}} \ \mathrm{d}V.
 \label{bilin_form}
\end{align}
and $((\cdot,\cdot))_{\manifold{M}, \mathbb{K}}\,:\, \Lambda^k\left(\manifold{M} \right ) \times \tilde{\Lambda}^{n-k}\left(\manifold{M} \right ) \rightarrow \mathbb{R}$ given by
\begin{align}
 (( {\mathbf{u}},{\mathbf{\tilde{v}}}))_{\manifold{M}, \mathbb{K}} := \int_{\manifold{M}}  \mathbf{u} \cdot \mathbb{K}^{-1}\mathbf{\tilde{v}} \ \mathrm{d}V.
\end{align}
For $\mathbf{q}\in \Lambda^2\left(\manifold{M} \right )$ and $p \in \tilde{\Lambda}^0\left(\manifold{M}\right )$ and homogeneous boundary values we have
\begin{equation}
\begin{aligned}
 \left( \left( \text{div} \mathbf{q}, \ {p} \right) \right)_{\manifold{M}} &=
 - \left( \left( \mathbf{q},\text{grad} p \right) \right)_{\manifold{M}} \\
 &= - \left( \left( \mathbf{q}, \mathbb{K}^{-1} \left[ \mathbb{K} \ \text{grad} p \right] \right) \right)_{\manifold{M}}  \\
 &=  \left( \left( \mathbf{q}, \text{grad}_{\mathbb{K}}^{*} \ p  \right)  \right)_{\manifold{M}, \mathbb{K}} .
\end{aligned}
\label{eq::adjointGradientDivergence}
\end{equation}
It is possible to define a new gradient operator,
\begin{align}
 \text{grad}_{\mathbb{K}}^{*} = - \mathbb{K} \ \text{grad}.
\label{def:grad_star}
\end{align}
\end{svgraybox}

\subparagraph{\textbf{Mixed formulation}}

Starting from \eqref{equation::darcyProblem} and making use of the bilinear maps defined above we have for all vectors $\boldsymbol{\tau}\in \Lambda^{n-1}\left(\manifold{M}\right )$ associated to outer-oriented surfaces
\begin{equation}
\begin{aligned}
\left( \left( \boldsymbol{\tau}, \vec{u} - \text{grad} p \right) \right)_{\manifold{M}} &= 0 \\
\Longleftrightarrow & \\
\left( \left( \boldsymbol{\tau}, \mathbb{K} \vec{u} - \mathbb{K}\text{grad} p \right) \right)_{\manifold{M},\mathbb{K}}  &= 0  \\
\stackrel{(\ref{def:grad_star})}{\Longleftrightarrow} & \\
\left( \left( \boldsymbol{\tau}, \mathbb{K} \vec{u} \right) \right)_{\manifold{M},\mathbb{K}} + \left( \boldsymbol{\tau} , \text{grad}^{*}_{\mathbb{K}} p \right)_{\manifold{M},\mathbb{K}} &= 0  \\
\stackrel{(1\mathrm{c}) \mbox{ and }(\ref{eq::adjointGradientDivergence})}{\Longleftrightarrow} & \\
\innerspace{\boldsymbol{\tau}}{\vec{q}}{\manifold{M},\mathbb{K}} + \left( \left( \text{div} \ \boldsymbol{\tau}, \ p \right) \right)_{\manifold{M}}  &= 0
\end{aligned}
\end{equation}
The constitutive equation is included in the last step by converting the bilinear form to a weighted inner product on $\Lambda^2\left(\manifold{M}\right )$ as defined in (\ref{def:weighted_inner}). For $(1\mathrm{b})$ we take the bilinear map for the divergence of a vector $\mathbf{q}$ associated with outer-oriented surfaces and an arbitrary scalar function defined in inner-oriented points, $\gamma \in \tilde{\Lambda}^0\left(\manifold{M} \right )$,
\begin{align}
\left ( \left ( {\text{div} \mathbf{q}},{\gamma} \right ) \right )_{\manifold{M}} &= \left ( \left ({\phi},{\gamma} \right ) \right )_{\manifold{M}}.
\end{align}

\begin{svgraybox}
The mixed formulation becomes: Find $\left(\mathbf{q}, p \right) \in \left\{ \Lambda^2\left(\manifold{M} \right) \times \tilde{\Lambda}^0\left(\manifold{M} \right) \right\}$, given $\phi \in \Lambda^3 \left(\manifold{M} \right)$, for all $\left( \boldsymbol{\tau}, \gamma \right) \in \left\{ \Lambda^2 \left(\manifold{M} \right) \times \tilde{\Lambda}^0 \left(\manifold{M} \right) \right\}$ such that,
\begin{align}
\inner{\tau}{\mathbf{q}}_{\manifold{M},\mathbb{K}} + \left( \left( \text{div} \ \boldsymbol{\tau}, \ p \right) \right)_{\manifold{M}}  &= 0 \label{sp1}\\
\left ( \left ( {\text{div} \mathbf{q}},{\gamma} \right ) \right )_{\manifold{M}} &= \left ( \left ({\phi},{\gamma} \right ) \right )_{\manifold{M}} \label{sp2}\;.
\end{align}
\end{svgraybox}

\subsection{Basis functions}

For the high order representation we use Lagrange, $l_{i} \left( \xi \right)$, and edge functions, $e_{i}\left( \xi \right)$. Lagrange polynomials interpolate nodal values. The edge functions, derived by Gerritsma \cite{gerritsma::edge_basis} are constructed such that when integrating over a line segment it gives one for the corresponding element and zero for any other line segment,
\begin{align}
l_{i} \left( \xi_{j} \right) = \delta_{i,j} \qquad
\int_{\xi_{j-1}}^{\xi_{j}} e_{i} \left( \xi \right) = \delta_{i,j}.
\end{align}
The relation between the Lagrange and the edge functions is given by,
\begin{align}
e_{i}\left( \xi \right) = \epsilon_{i} \left( \xi \right) d\xi, \quad \text{with} \quad \epsilon_{i}\left( \xi \right) = - \sum_{k=0}^{i-1} \frac{d l_{k}}{d \xi}.
\end{align}
Note that this definition implies
\begin{equation} \frac{d l_{i}}{d \xi} = e_i(\xi) - e_{i+1}(\xi) \;.\label{relation_deriv_edge}\end{equation}

Extension to the multidimensional is obtained by means of tensor products. For more details see \cite{kreeft2011mimetic}.

\subsection{Mimetic discretization in 2D}
\label{Section::discreteIntegrals}
\subparagraph{\textbf{Expansion of unknowns in $\realNumber^{2}$}}
Let $\mathbf{q} \in \Lambda^1\left(\manifold{M} \right)$ be expanded as,
\begin{align}
 \mathbf{q}_h = \left[ \begin{array}{c}
                   \sum\limits_{i=0}^{N}\sum\limits_{j=1}^{N} q_{i,j}^{x} l_{i}\left( \xi \right) e_{j} \left( \eta \right) \\
                   \sum\limits_{i=1}^{N}\sum\limits_{j=0}^{N} q_{i,j}^{y} e_{i}\left( \xi \right) l_{j} \left( \eta \right)
                  \end{array}
 \right]\;,
\end{align}
and the pressure, $p_h \in \tilde{\Lambda}^0\left(\manifold{M} \right)$ as,
\begin{align}
 p_h = \sum\limits_{i=1}^{N}\sum\limits_{j=1}^{N} p_{i,j} \epsilon_{i}\left( \xi \right) \epsilon_{j} \left( \eta \right).
\end{align}

\begin{figure}[ht]
\centering
\subfigure{
\includegraphics[width=1\textwidth]{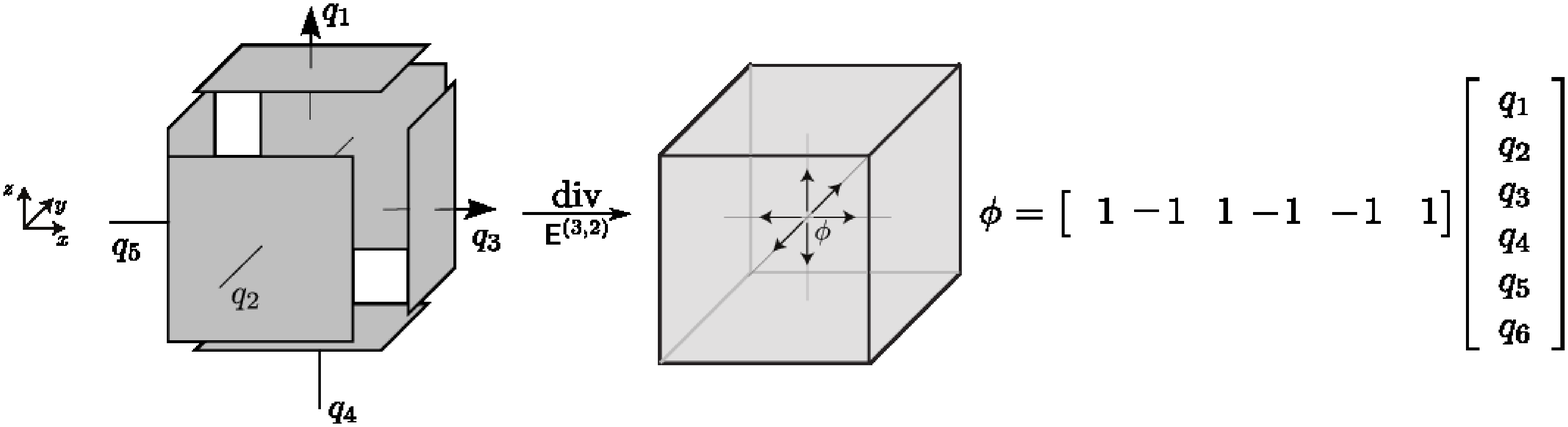}
}
\subfigure{
\includegraphics[width=0.85\textwidth]{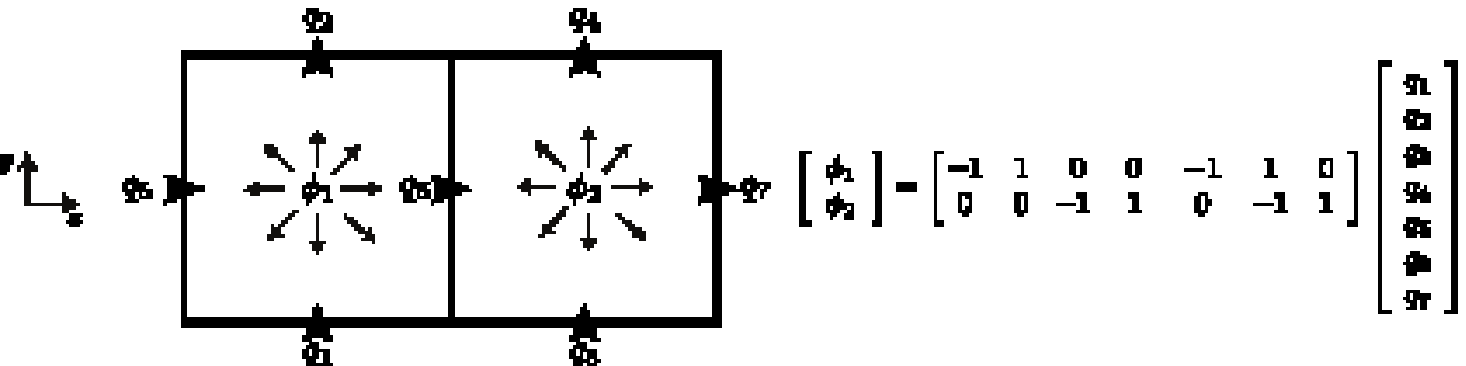}
}
\caption{Discrete representation of the action of the divergence in $\realNumber^{3}$ and $\realNumber^{2}$.}
\label{fig:discreteDivergence}
\end{figure}

\subparagraph{\textbf{Discrete divergence in $\realNumber^{2}$}}
The divergence of $\mathbf{q}_h$ is then given by
\begin{align}
 \text{div} \ \vec{q}_h &= \sum\limits_{i=1}^{N}\sum\limits_{j=1}^{N} \left( q_{i,j}^{x} - q_{i-1,j}^{x} + q_{i,j}^{y} - q_{i,j-1}^{y} \right) e_{i}\left( \xi \right) e_{j} \left( \eta \right) \;,
 \label{discrete_div}
\end{align}
where we repeatedly used (\ref{relation_deriv_edge}). The scalar $\phi \in \Lambda^2\left(\manifold{M} \right)$ associated with outer-oriented volumes is expanded as
\begin{align}
 \phi_h &= \sum\limits_{i=1}^{N}\sum\limits_{j=1}^{N} \phi_{i,j} e_{i}\left( \xi \right) e_{j} \left( \eta \right) \;.
 \label{discrete_phi}
\end{align}
Equating (\ref{discrete_div}) and (\ref{discrete_phi}) yields
\begin{align}
 \sum\limits_{i=1}^{N}\sum\limits_{j=1}^{N} \phi_{i,j} \bcancel{ e_{i}\left( \xi \right) e_{j} \left( \eta \right) } &= \sum\limits_{i=1}^{N}\sum\limits_{j=1}^{N} \left( q_{i,j}^{x} - q_{i-1,j}^{x} + q_{i,j}^{y} - q_{i,j-1}^{y} \right) \bcancel{ e_{i}\left( \xi \right) e_{j} \left( \eta \right)} \label{eq::basiscancel} \\
 \left[\phi \right] &= \incidencederivative{2}{1} \left[\vec{q} \right]. \label{topological_relation}
\end{align}
We see that the basis functions cancel from this relation. The matrix $\incidencederivative{2}{1}$ relates the fluxes $q_{i,j}^x$ and $q_{i,j}^y$ to the volume integral $\phi_{i,j}$, as depicted in Figure~\ref{fig:discreteDivergence}. This fully discrete equation is a restatement of the integral relation (\ref{integral_divergence}). The matrix $\incidencederivative{2}{1}$ only contains the values $-1$, $0$ and $1$ and is fully determined by the grid, see \cite{kreeft2011mimetic}. This is an incidence matrix showing the topological nature of the discrete divergence.

\begin{svgraybox}
If we insert the expansions of our unknowns in (\ref{sp1}) and (\ref{sp2}) we obtain in $\mathbb{R}^n$ the saddle point problem given by,
\begin{align}
\left[
 \begin{array}{cc}
  × \matrixoperator{M}^{(n-1)}_{ \mathbb{K} } & \left(\incidencederivative{n}{n-1}\right)^{T} \matrixoperator{M}^{(n)} \\
  × \matrixoperator{M}^{(n)} \incidencederivative{n}{n-1} & 0
 \end{array}
\right]
\left[
 \begin{array}{cc}
  × \vec{q}\\
  × p
 \end{array}
\right] =
\left[
 \begin{array}{cc}
  × 0\\
  × \matrixoperator{M}^{(n)} \phi
 \end{array}
\right],
\label{eq::systemMatrixDarcy}
\end{align}
where $\matrixoperator{M}^{(k)}$ is the symmetric mass matrix obtain from the bilinear pairing between variables associated with outer-orientation and inner-orientation, (\ref{bilin_form}), $\matrixoperator{M}^{(n-1)}_{ \mathbb{K} }$ is the mass matrix obtained from the weighted inner product (\ref{def:weighted_inner}) and $\incidencederivative{n}{n-1}$ the incidence matrix which relates fluxes over surfaces to volumes. The resulting system \eqref{eq::systemMatrixDarcy} is symmetric.
\end{svgraybox}

%

The pressure which is represented on an inner-oriented grid (which is not explicitly constructed in this single grid approach) is pre-multiplied by $\matrixoperator{M}^{(n)}$ to represent it on the outer-oriented grid.

\section{NUMERICAL RESULTS}
\label{Section::Numerical}

The method derived in this paper respects the geometric nature of the problem. However, it is crucial to verify the numerical benefits of this approach. This section presents $hp-$convergence studies for anisotropic permeability.

\subsection{Manufactured solution - Anisotropic permeability}

\begin{figure}[ht]
\centering
\subfigure{
\includegraphics[width=0.4\textwidth]{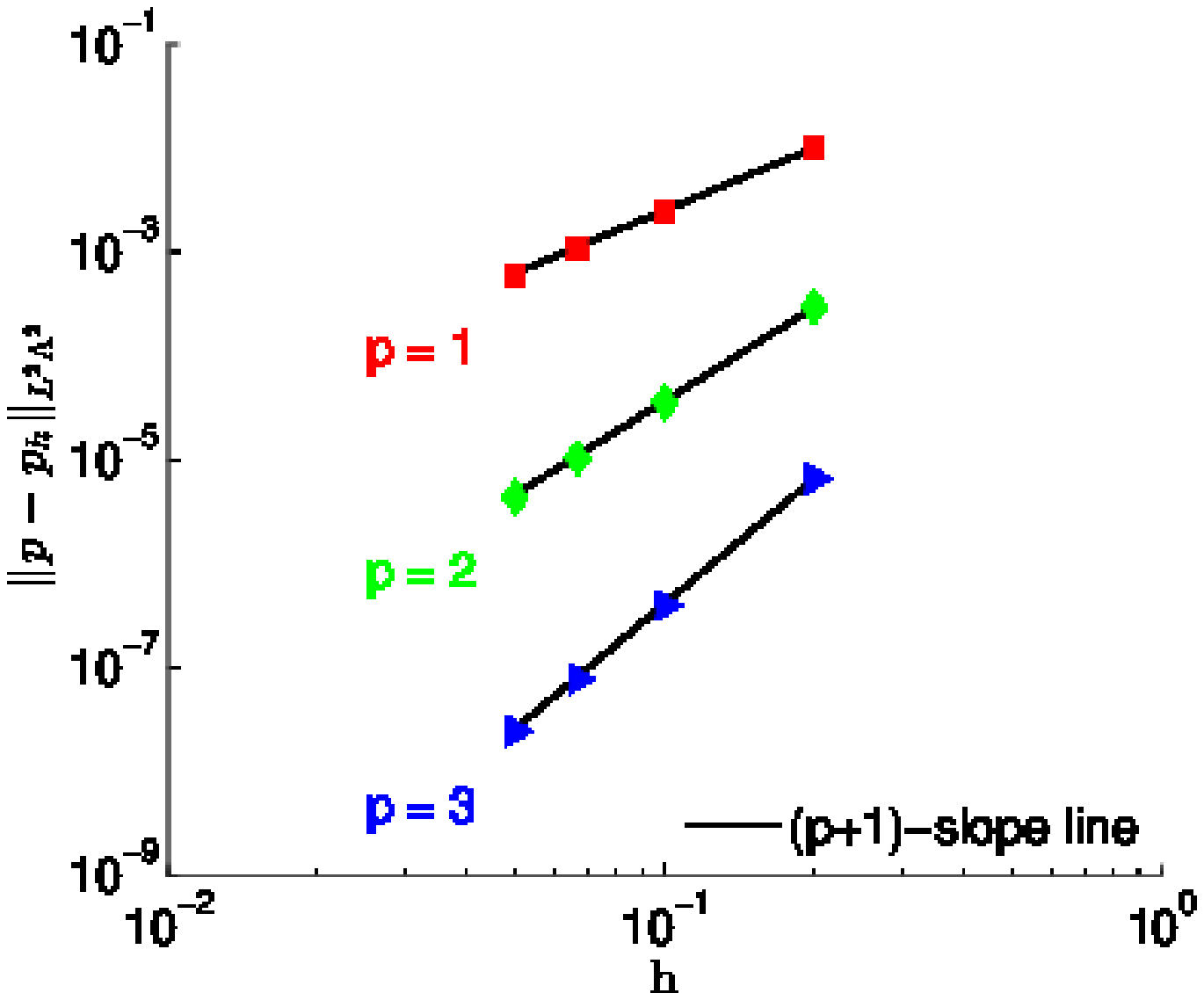}
\label{fig:}
}
\subfigure{
\includegraphics[width=0.4\textwidth]{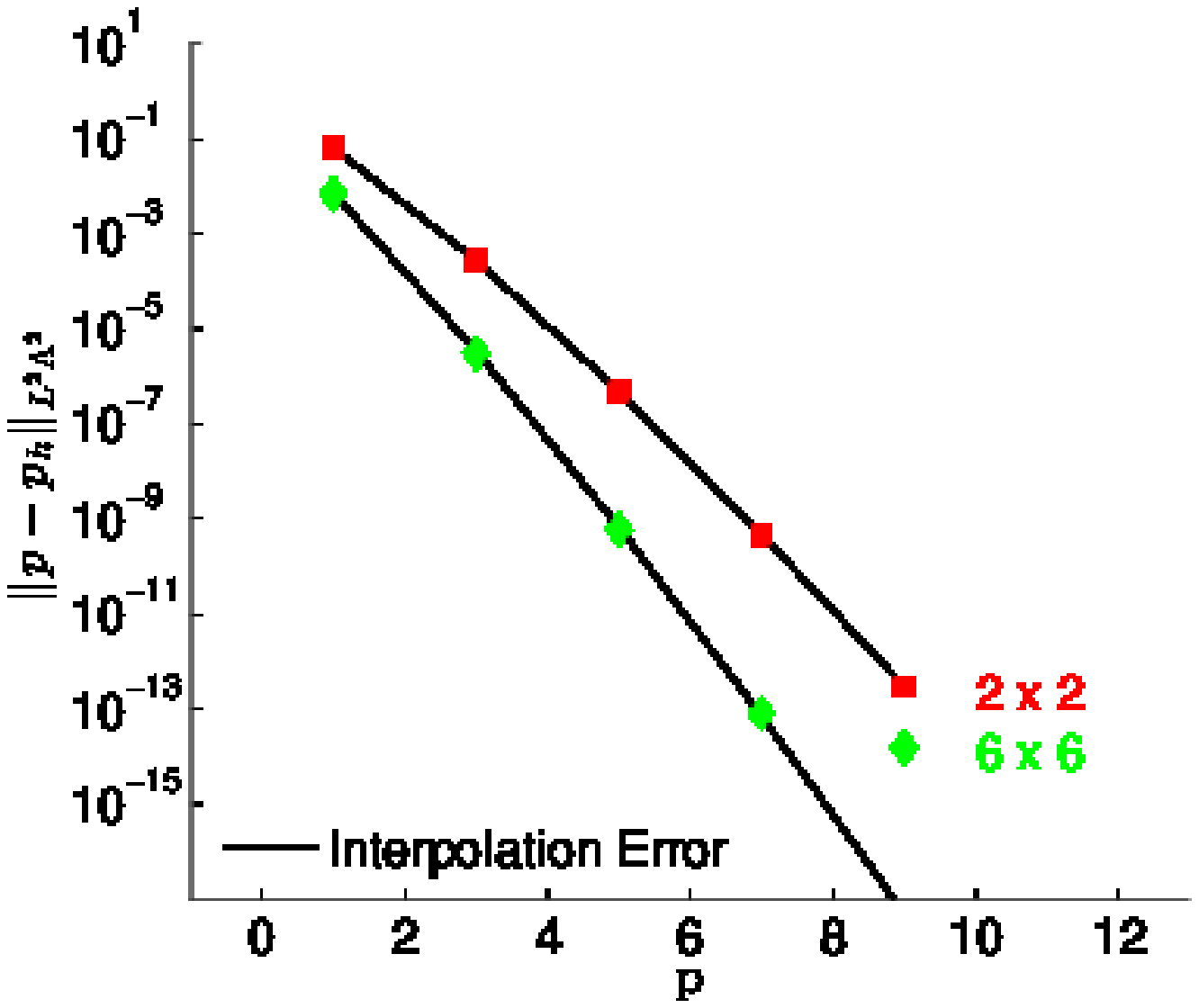}
\label{fig:}
}
\caption{Plots of the $h-$ and $p-$convergence for anisotropic permeability given in \eqref{equation::AnisotropicPermeability}.}
\label{fig:hpConvergencePressurePatch}
\end{figure}

The first test case assesses the convergence for $h-$ and $p-$refinement of the mixed mimetic spectral element method applied to the Darcy model. This is a benchmark problem presented in \cite{hyman1997numerical}. The problem is defined on a unit square, $\Omega = \left[-1,1 \right]^{2}$, with Cartesian coordinates with permeability given by,
\begin{equation}
\mathbb{K} = \begin{bmatrix}
               2 & 1 \\
               1 & 2
              \end{bmatrix}
\label{equation::AnisotropicPermeability}
\end{equation}
and the right hand side, $\phi \in L^{2}\left({\manifold{M}}\right)$ given by,
\begin{align}
 \kdifform{\phi}{2} = 2\left( 1 + x^2 + xy + y^2 \right)e^{xy} \ \mathrm{d}x \mathrm{d}y.
\end{align}
This results in an exact solution for pressure $p \in \tilde{\Lambda}^0\left(\manifold{M} \right)$ given by,
\begin{align}
 \kdifform{p}{0} = e^{xy}
\end{align}
\figref{fig:hpConvergencePressurePatch} shows the $h-$ and $p-$convergence for the pressure in straight mesh. For the  $h-$convergence the expected rate of convergence is of $(p+1)$, where $p$ is the polynomial degree. The solid line for the interpolation error in the $p$-convergence plot is the $L^{2}$-error from interpolating the exact solution, the solution converges exponentially. Both the numerical solution and the interpolated exact solution converge exponentially.

\subsection{Layered medium}

\begin{figure}[ht]
\centering
\subfigure{
\includegraphics[width=0.4\textwidth]{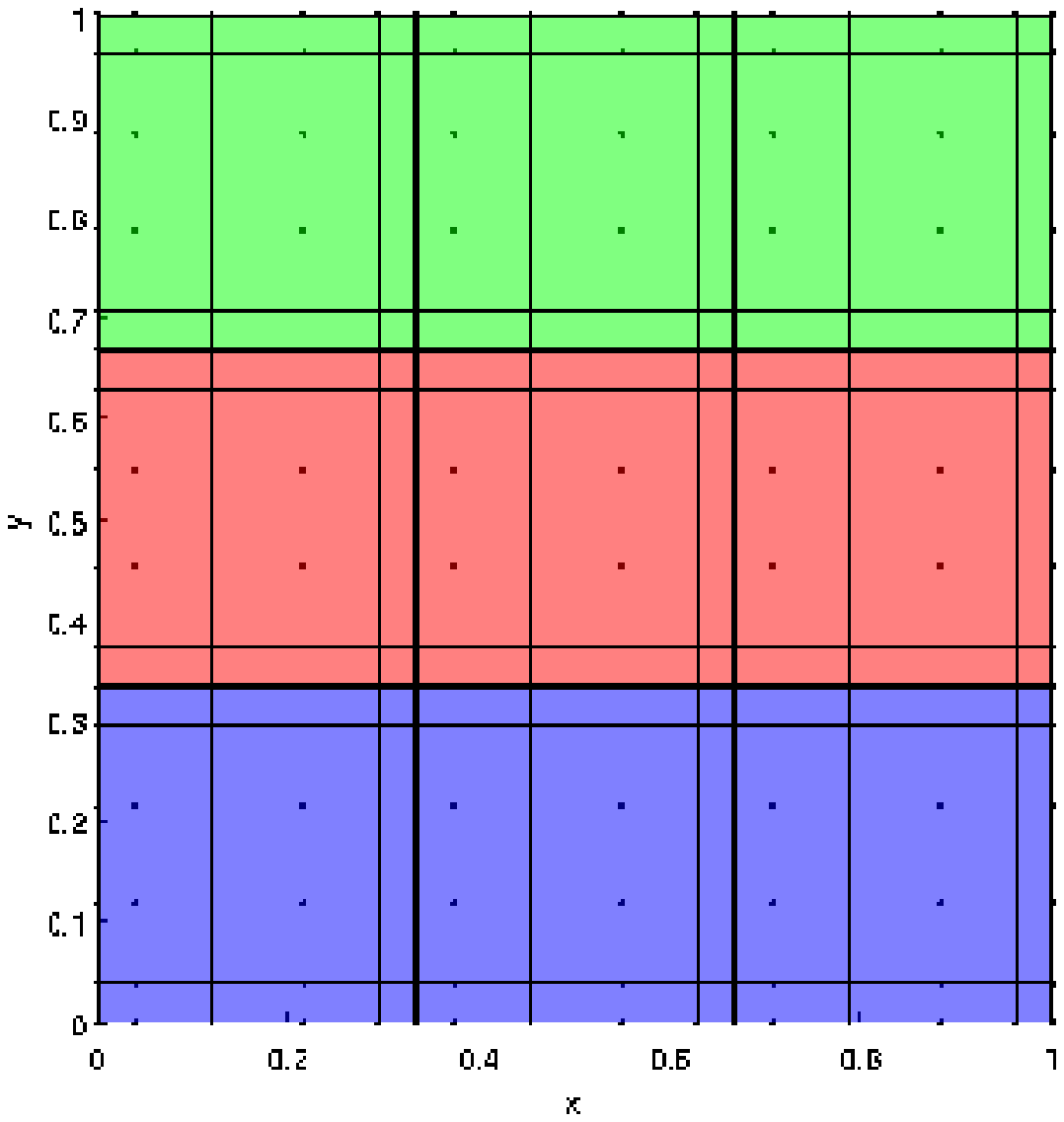}
\label{fig:crazyMesh_cc00}
}
\subfigure{
\includegraphics[width=0.4\textwidth]{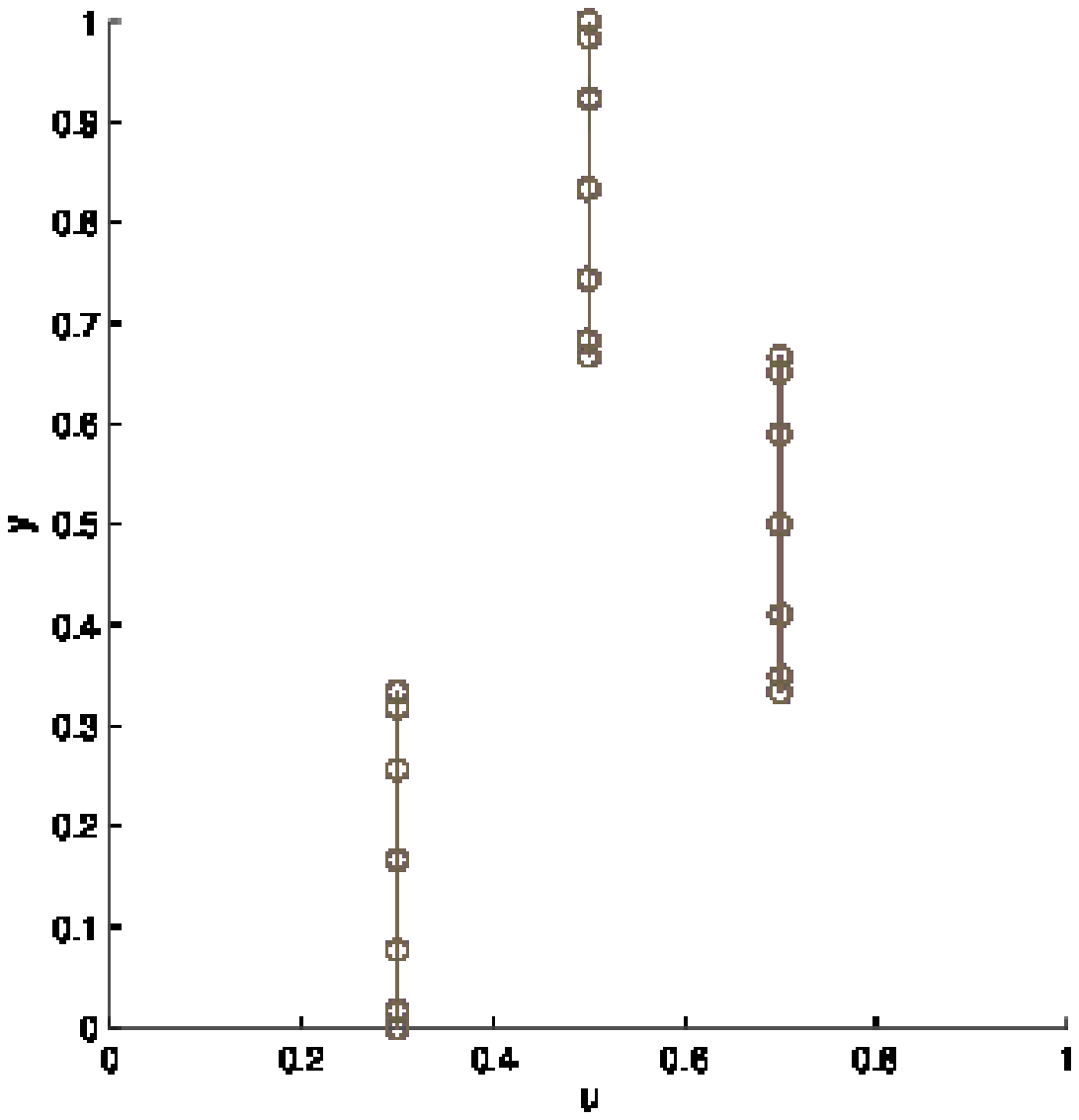}
\label{fig:crazyMesh_cc025}
}
\caption{Layered medium with 3 different permeabilities, $\alpha$ ranges top to bottom from $0.5$, $0.7$ and $0.3$. The left figure shows a velocity profile at an arbitrary $x$. Note the numerical solution in the Gauss-Lobatto nodes.}
\label{fig:layeredMedium}
\end{figure}

A classical benchmark for Darcy flow codes is the piecewise constant permeability in a square \cite{masud2002stabilized}. Such a medium is called \textit{layered} medium.
\begin{equation}
\mathbb{K} = \begin{bmatrix}
               \alpha & 0 \\
               0 & \alpha
              \end{bmatrix} \quad \alpha = \begin{cases}
                                            0.3 & \text{if $y \leq \frac{1}{3}$}\\
                                            0.7 & \text{if $\frac{1}{3} < y \leq \frac{2}{3}$}\\
                                            0.5 & \text{if $y > \frac{2}{3}$}
                                           \end{cases}
\end{equation}
The fluid comes into the domain from the left to the right. Since the pressure depends linearly on $x$, horizontal constant velocity is expected in each layer, \figref{fig:layeredMedium}.

\begin{acknowledgement}
The authors gratefully acknowledge the funding received by FCT - Foundation for science and technology Portugal through SRF/BD/36093/2007 and SFRH/BD/79866/2011 and the anonymous reviewers for their helpful comments.
\end{acknowledgement}

\bibliographystyle{plain}


\begin{thebibliography}{10}

\bibitem{arnold2010finite}
D.~Arnold, R.~Falk, and R.~Winther.
\newblock {Finite element exterior calculus: from Hodge theory to numerical
  stability}.
\newblock {\em American Mathematical Society}, 47(2):281--354, 2010.

\bibitem{jerome2012analysis}
J.~Bonelle and A.~Ern.
\newblock Analysis of compatible discrete operator schemes for elliptic
  problems on polyhedral meshes.
\newblock {\em arXiv preprint arXiv:1211.3354}, 2012.

\bibitem{bossavit2005discretization}
A~Bossavit.
\newblock Discretization of electromagnetic problems.
\newblock {\em Handbook of numerical analysis}, 13:105--197, 2005.

\bibitem{BrezziBuffaLipnikov2009}
F~Brezzi, A~Buffa, and K~Lipnikov.
\newblock {{M}imetic finite differences for elliptic problems}.
\newblock {\em Mathematical Modelling and Numerical Analysis}, 43(2):277--296,
  2009.

\bibitem{burke1985applied}
W.~L. Burke.
\newblock {\em {Applied differential geometry}}.
\newblock Cambridge Univ Pr, 1985.

\bibitem{desbrun2005discrete}
M.~Desbrun, A.~Hirani, M.~Leok, and J.~Marsden.
\newblock {Discrete exterior calculus}.
\newblock {\em Arxiv preprint math/0508341}, 2005.

\bibitem{frankel}
T.~Frankel.
\newblock {\em {{T}he {G}eometry of {P}hysics}}.
\newblock Cambridge University Press, 2nd edition, 2004.

\bibitem{gerritsma::edge_basis}
M~Gerritsma.
\newblock Edge functions for spectral element methods.
\newblock {\em Spectral and High Order Methods for Partial differential
  equations, Eds J.S. Hesthaven \& E.M. R{\o}nquist, Lecture Notes in
  Computational Science and Engineering}, 76.

\bibitem{gerritsmaICOSAHOM}
M.~Gerritsma, R.~Hiemstra, J.~Kreeft, A.~Palha, P.~Pinto~Rebelo, and
  D.~Toshniwal.
\newblock The geometric basis of numerical methods.
\newblock {\em Proceedings ICOSAHOM 2012 (this issue)}, 2012.

\bibitem{herbin2008benchmark}
R.~Herbin and F.~Hubert.
\newblock Benchmark on discretization schemes for anisotropic diffusion
  problems on general grids.
\newblock {\em Finite volumes for complex applications V}, pages 659--692,
  2008.

\bibitem{HiemstraICOSAHOM}
R.~Hiemstra and M.~Gerritsma.
\newblock High order methods with exact conservation properties.
\newblock {\em Proceedings ICOSAHOM 2012 (this issue)}, 2012.

\bibitem{hiemstra2012high}
R.~Hiemstra, R.~Huijsmans, and M.~Gerritsma.
\newblock High order gradient, curl and divergence conforming spaces, with an
  application to compatible isogeometric analysis.
\newblock {\em Submitted to J. Comp Phys., arXiv preprint arXiv:1209.1793},
  2012.

\bibitem{Hirani_phd_2003}
A.~Hirani.
\newblock {\em {{D}iscrete {E}xterior {C}alculus}}.
\newblock PhD thesis, California Institute of Technology, 2003.

\bibitem{hirani2008numerical}
A.~Hirani, K.~Nakshatrala, and J.~Chaudhry.
\newblock Numerical method for {D}arcy flow derived using discrete exterior
  calculus.
\newblock {\em arXiv preprint arXiv:0810.3434}, 2008.

\bibitem{hyman1997numerical}
J.~Hyman, M.~Shashkov, and S.~Steinberg.
\newblock The numerical solution of diffusion problems in strongly
  heterogeneous non-isotropic materials.
\newblock {\em Journal of Computational Physics}, 132(1):130--148, 1997.

\bibitem{KreeftICOSAHOM}
J.~Kreeft and M.~Gerritsma.
\newblock Higher-order compatible discretization on hexahedrals.
\newblock {\em Proceedings ICOSAHOM 2012 (this issue)}, 2012.

\bibitem{kreeft::stokes}
J.~Kreeft and M.~Gerritsma.
\newblock Mixed mimetic spectral element method for stokes flow: a pointwise
  divergence-free solution.
\newblock {\em Journal of Computational Physics}, 2012.

\bibitem{kreeft2012priori}
J.~Kreeft and M.~Gerritsma.
\newblock A priori error estimates for compatible spectral discretization of
  the stokes problem for all admissible boundary conditions.
\newblock {\em arXiv preprint arXiv:1206.2812}, 2012.

\bibitem{kreeft2011mimetic}
J.~Kreeft, A.~Palha, and M.~Gerritsma.
\newblock {Mimetic framework on curvilinear quadrilaterals of arbitrary order}.
\newblock {\em Submitted to FoCM, Arxiv preprint arXiv:1111.4304}, 2011.

\bibitem{masud2002stabilized}
A.~Masud and T.J.R. Hughes.
\newblock A stabilized mixed finite element method for darcy flow.
\newblock {\em Computer Methods in Applied Mechanics and Engineering},
  191(39):4341--4370, 2002.

\bibitem{palhaAdvection}
A.~Palha, P.~Pinto~Rebelo, and M.~Gerritsma.
\newblock Mimetic spectral element solution for conservative advection.
\newblock {\em Proceedings ICOSAHOM 2012 (this issue)}, 2012.

\bibitem{palhaLaplaceDualGrid}
A.~Palha, P.~Pinto~Rebelo, R.~Hiemstra, J.~Kreeft, and M.~Gerritsma.
\newblock Physics-compatible discretization techniques on single and dual
  grids, with application to the {P}oisson equation of volume forms.
\newblock {\em Submitted to J. Comp Phys.}, 2012.

\bibitem{tonti1975formal}
E~Tonti.
\newblock {On the formal structure of physical theories}.
\newblock {\em preprint of the Italian National Research Council}, 1975.

\bibitem{ToshniwalICOSAHOM}
D.~Toshniwal, R.H.M. Huijsmans, and M.~Gerritsma.
\newblock A geometric approach towards momentum conservation.
\newblock {\em Proceedings ICOSAHOM 2012 (this issue)}, 2012.

\end{thebibliography}

\end{document}